\documentclass[12pt]{article}

\usepackage{amsmath,amssymb}

\def\CC{{\mathbb C}}
\def\ZZ{{\mathbb Z}}

\def\RR{{\mathbb R}}
\def\QQ{{\mathbb Q}}
\def\DD{{\mathbb D}}
\def\PP{{\mathbb P}}

\def\ZZ{{\mathbb Z}}

\newcommand {\BOX} {\rule{2mm}{2mm}}

\newcommand {\be} {\begin{equation}}
\newcommand {\eeqn} {\end{equation}}
\newcommand {\bea} { \begin{eqnarray}}
\newcommand {\eea} {\end{eqnarray}}
\newcommand {\beas} { \begin{eqnarray*}}
\newcommand {\eeas} {\end{eqnarray*}}
\newcommand {\ra} {\rightarrow}

\newtheorem {lemma} {LEMMA} [section]
\newtheorem {theorem}[lemma]{THEOREM} 
\newtheorem {prop}[lemma]{PROPOSITION}
\newtheorem {cor}[lemma]{COROLLARY}
\newtheorem {definition}[lemma]{DEFINITION}
\newtheorem {question}[lemma]{QUESTION}

\numberwithin{equation}{section}

\begin{document}


\title{Double sections, dominating maps\\  
and the Jacobian fibration}
\author{Gregery T. Buzzard\thanks{Partially supported by an NSF grant.} $\ $
and Steven S. Y. Lu\thanks{Partially supported by an NSERC grant.}}
\date{}
\maketitle


\section{Introduction}

As is well-known, there exist nonconstant holomorphic maps from the
plane into the Riemann sphere $\PP^1$ minus two points,
the simplest example of which is an explicit realization of
the uniformization map given by applying the exponential map
and then composing with a Mobius transformation taking $0$ and $\infty$
to the two given punctures.  Likewise, we can map the plane into the
sphere minus one point by simply applying directly
a Mobius transformation taking $\infty$ to this puncture.  

In this paper we prove two parametrized versions of this uniformization
result using different methods.  We first present a slightly 
easier version using
complex analysis, which allows us to construct the uniformizing maps
more or less explicitly in terms of the initial coordinates given. 
Then we give the more general and coordinate invariant version, 
which we prove by extending Kodaira's theory
of the Jacobian fibration to a family of singular curves constructed
via algebraic geometry.

Finally, using the results obtained with the Jacobian fibration, we obtain
two equivalent conditions for an complex surface nonhyperbolically
fibered over a complex curve to be holomorphically dominable by $\CC^2$.  In
theorem~\ref{thm:surf} we show that this dominability is equivalent to
the apparently weaker condition of the existence of a Zariski dense
image of $\CC$ and equivalent to the quasiprojectivity of the base
curve together with the nonnegativity of the orbifold Euler
characteristic.  

We start with the following definition.

\begin{definition}
Let $R$ be a Riemann surface (also called a complex curve) and let $D
\subset R \times \PP^1$ be a 1-dimensional complex subvariety
such that for generic $z \in R$ the fiber $D_z = D \cap (\{z\} \times
\PP^1)$ has two points.
Then $D$ is called a double section over $R$.
\end{definition}

Our first parametrized version of uniformization is that, given a double
section $D$ over a noncompact Riemann surface $R$, 
there is a fiber preserving holomorphic map 
from $R\times\CC$ to the complement of $D$ which on each fiber is a
uniformizing map as described above.

\begin{theorem} \label{thm:dominate}
Let $D$ be a double section over a noncompact Riemann surface $R$.  
Then there exists a holomorphic
map $F: R\times\CC \rightarrow (R \times \PP^1) \setminus D$ of the form
$F(z,w) = (z, H(z,w))$.  Moreover, if $\#D_z = 2$, then $H(z, w) =
M_z(\exp(c_z w))$, where $M_z$ is a Mobius transformation mapping 0 and
$\infty$ to the two points in $D_z$ and $c_z$ is a nonzero constant
depending on $z$, while if $\#D_z = 1$, then $H(z,w)
= M_z(w)$ where $M_z$ is a Mobius transformation mapping $\infty$ to
the point in $D_z$.
\end{theorem}

We will obtain this theorem by first proving a weaker result.

\begin{theorem} \label{thm:section}
Let $D$ be a double section over a noncompact Riemann surface $R$.  
Then there exists a holomorphic
function $\sigma: R \rightarrow \PP^1$ whose graph, a
holomorphic section over $R$, is disjoint from $D$.
\end{theorem}

As will be shown in section 5, if we replace $D$ by a triple section, 
then this theorem is no longer true.
\smallskip

Note that any double section $D$ necessarily has the form
$$D = \{(z,w): a(z)w^2 + b(z) w + c(z) = 0\}$$ for holomorphic $a, b, c$.
But, as a corollary of theorem~\ref{thm:section}, we can always 
put $D$ into a normal form in the case $R$ is noncompact:

\begin{cor}  \label{cor:normal}
Let $D$ be a double section over a noncompact Riemann surface $R$.  
Then after a biholomorphic
change of variables on $R \times \PP^1$ preserving each fiber, 
there is a holomorphic function $h$ on $R$ such that $D$ has the form
$$ D = \{(z,w): w^2 = h(z)\}. $$
\end{cor}

We now give the more general version of our simultaneous 
uniformization result.  Let $p: P\rightarrow R$ 
define a $\PP^1$ bundle over a Riemann surface $R$. We may define a
double section $D$ of $p$ as before. However, we will still reserve
the term ``double section over $R$'' to mean that $p$ defines a trivial
holomorphic $\PP^1$ bundle with an explicit trivialization, which can 
always be realized when $R$ is noncompact. Associated with $(P,D)$ is a
Jacobian fibration, to be introduced in section 5,
which is a complex analytic surface
obtained as a holomorphic quotient of a line bundle
over $R$. The following theorem states that if $p$ has
a section avoiding $D$, which as seen above is true when $R$ is noncompact,
then $P\setminus D$ is isomorphic to this Jacobian fibration. 

\begin{theorem}\label{thm:Jac}
Let $p: P\rightarrow R$
define a holomorphic $\PP^1$ bundle over a Riemann surface $R$ with a 
double section $D$. If $p$ has a section $\sigma$ that avoids $D$
(which is necessarily the case when $R$ is noncompact), then
there is a holomorphic line bundle $L$ over $R$ and a holomorphic
map $F: L\rightarrow P\setminus D$ preserving each fiber such that 
$F_z:L_z\rightarrow \PP^1\setminus D_z$ is isomorphic to the universal covering
map for all $z\in R$. Moreover, $F$ and $L$ are canonically determined by 
$(p,D,\sigma)$ and $F$ maps the zero-section of $L$ to $\sigma$.
\end{theorem}

Readers interested in the more abstract methods of algebraic geometry
and the precise setting as detailed in the proof of this theorem
should see section 3. However, section 3 can be skipped without losing
continuity. This version implies the previous theorems
by using this theorem in a local way, then solving an additive 
Cousin problem to obtain a global section as in theorem 1.3.  
This is done in section 3 as a corollary to proposition~\ref{prop:L3}.

	The Jacobian fibration to be introduced in 
section 3 is analogous to and indeed based on the one
introduced by Kodaira \cite{Kod,Kod1,Kod2} in his study of elliptic
fibrations.  This work of Kodaira served as an important foundation for the
modern classification theory of compact complex surfaces.  The
Jacobian fibration introduced in the current paper is intended to
serve as a precursor to a wider theory in which noncompact
compactifiable complex manifolds may be studied more effectively and in
greater detail and to provide motivation for studying the Jacobian
fibration in the realm of families of singular varieties.

The motivation for this paper came from ideas contained in \cite{bl}.
In that paper we showed that for many algebraic surfaces, the
existence of a Zariski dense image of $\CC$ is equivalent to the
existence of a dominating map from $\CC^2$ into the surface (a
dominating map is a holomorphic map with full rank generically).
However, in that paper we did not resolve completely the question of
the existence of a Zariski dense image of $\CC$ in the complement of a
general double section.  The theorems above resolve this problem
completely. 

Additionally, in keeping with the general
strategy of \cite{bl} to deal with the dominability question
for arbitrary complex surfaces, 
we obtain the following result. It states
roughly that given any discrete set of jet prescriptions of sections
over $R$, we can construct a dominating map, $F$, whose image avoids the 
double section outside the given discrete set 
while at the same time realizing the jet prescriptions
in the following sense:  if $c$ is a point in $R$ with a corresponding
jet prescription and $\sigma$ is a holomorphic function in a
neighborhood of $c$, then $F (z, \sigma(z))$ has the given jet
prescription at $z=c$. 

\begin{theorem}   \label{thm:jet}
Let $D$ be a double section over a noncompact Riemann surface $R$.  
Let $E$ be a discrete subset of
$R$, and for each $c \in E$, let $P_c(z)$ be a Laurent polynomial in
powers of $z-c$.  Then there exists a holomorphic map $F$ from $R\times \CC$
into $R \times \PP^1$ of the form $F(z,w) = (z, H(z,w))$ such that if
$c \in E$ and $\sigma$ is a holomorphic function in a neighborhood of
$c$, then the Laurent expansion of $F(z, \sigma(z))$ agrees with
$P_c(z)$ up to order $\deg(P_c)$, while if $c \notin E$, then $F(c,
\cdot)$ maps $\CC$ onto $\PP^1\setminus D_c$ as in
theorem~\ref{thm:dominate}.  
\end{theorem}

	Coupled with results in \cite{bl} and some geometry 
in section 4.1 on orbifold Riemann surfaces, 
this allows us to give a complete geometric classification of 
nonhyperbolically fibered complex surfaces which are dominable 
by $\CC^2$, i.e, which admit holomorphic maps from $\CC^2$ having full
rank somewhere. To state it, recall that a fibered surface $X$ 
is given by a surjective holomorphic map $f$ 
from the complex surface $X$ to a complex curve $C$ that extends to
a proper holomorphic map on a complex surface $\bar X$ containing $X$ 
as a Zariski open subset. 
This gives a set $M_s$ of multiplicities for
each point $s$ of $C$, comprising the various degrees of 
vanishing of $f$ (or multiplicities of $f$) at points of $X_s=f^{-1}(s)$
when $s$ is given as the origin in a local coordinate on $C$.
Recall also that the classical multiplicity $n_s$  of the fiber $X_s$
is given as $\gcd M_s$. We define the new multiplicity $m_s$ at $s$
to be $\min M_s$. Then our result says that if
the general fiber of $f$ is nonhyperbolic, then
$X$ is dominable by $\CC^2$ if and only if
$C$ together with the $\QQ$-divisor $A=\sum (1-1/m_s)s$
is not hyperbolic as a marked Riemann surface. The latter happens precisely
when $C$ is quasiprojective (with smooth compactification $\bar C$)
and $\chi(C,A)\geq 0$, where the orbifold Euler characteristic is given by
$$
\chi(C,A):= 2-2g(\bar C)-\#(\bar C\setminus C) - \deg A = e(C) - \deg A
$$
as in \cite{bl}.

\begin{theorem}\label{thm:surf}
Let $f:X\rightarrow C$ be a surjective holomorphic map defining $X$
into a fibered surface over the 
curve $C$. This imposes a branched orbifold structure on
$C$ specified by a $\QQ$-divisor $A=\sum (1-1/m_s)s$, where $m_s$ is
the minimum multiplicity of $f$ at points above $s$, for each point $s\in C$. 
Assume that the general fiber 
$X_c$ is not hyperbolic, or equivalently, that $e(X_c)\geq 0$.
Then the following are equivalent.

\begin{itemize}
\item[\rm (a)] $X$ is dominable by $\CC ^ {2}$.
\item[\rm (b)]  $C$ is quasiprojective and  $\chi(C,A)\geq 0$.
\item[\rm (c)]  There exists a holomorphic map of $\CC$ to $X$ whose
image is Zariski dense. 
\end{itemize}
\end{theorem}

This paper is organized as follows. Section 2 proves 
theorem~\ref{thm:dominate} and theorem~\ref{thm:jet} by first
establishing theorem~\ref{thm:section} in section 2.1.  From this,
theorems~\ref{thm:dominate} and theorem~\ref{thm:jet} and
corollary~\ref{cor:normal} follow naturally, as seen in sections 2.2 and 2.3.
The methods presented in section 2
are completely concrete and allow us
to write the dominating maps explicitly in terms of the initial coordinates.
Section 3 presents the more abstract setting of Jacobian fibrations
associated with our problem and uses this naturally coordinate
invariant setting to prove theorem~\ref{thm:Jac}.  Section 4 
recalls the setting for theorem~\ref{thm:surf} from \cite{bl},
establishes a generalized orbifold version of Ahlfors' lemma, and 
uses theorem~\ref{thm:Jac} and the $n$th-root 
fibration to establish theorem~\ref{thm:surf}.
Finally, section 5 provides some further remarks, including why 
our result is sharp to some extent, and gives further directions for
investigation. 

We are grateful to Bernard Shiffman for his aid in establishing this 
fruitful collaboration, and in particular for posing the problem that
started it. 


\section{Dominating the complement of a double section}
Throughout this section, $R$ will be assumed to be a noncompact 
Riemann surface, thus allowing us to use basic interpolation theorems 
on $R$ to construct the necessary maps.

\subsection{A holomorphic section}

In this section we prove theorem~\ref{thm:section}.  Let $D$ be a
double section in $R \times \PP^1$.  Since there are only countably
many $z \in R$ for which $\#D_z = 1$, we can apply a map of the form
$(z,w) \mapsto (z, M(w))$ for some Mobius transformation $M$ to ensure
that if $\#D_z = 1$, then $\infty \notin D_z$. 

Let $E_\infty$ be the
set of $z \in R$ such that $\infty \in D_z$.  Then $E_\infty$ is
discrete, so we may apply a map of the form $(z, w) \mapsto (z,
w-c)$, where $c$ is a constant such that if $z \in
E_\infty$ and $w \in D_z$, then $w-c \neq 0$.  Hence we may assume
that $D$ has been normalized so that if $z_0 \in
E_\infty$, then there is a neighborhood of $z_0$ such that in this
neighborhood, $D$ is the union of the graphs of two functions of the
form $u_0(z)$ and $u_1(z)/(z-z_0)^m$, where $u_j(z_0) \neq 0$ and $m
\geq 1$. 

Next, by \cite[theorems 26.3, 26.5]{forster}, the
classical theorems of Weierstrass and Mittag-Leffler are valid on
$R$, so there exists $f_1$ holomorphic on $R$ such that if $z_0 \in E_\infty$
with $u_0, u_1$ as described above, then $f_1(z) = (z-z_0)^m +
O(|z-z_0|^{m+1})$.  Let $F_1(z,w) = (z, f_1(z) w)$, and let $C =
F_1(D)$.  Then $C$ is a double section and is contained in $R \times
\CC$.  Moreover, if $z_0 \in E_\infty$, then $C$ is locally the union
of the graphs of $v_0(z) = u_0(z)(z-z_0)^m$ and $v_1(z) = u_1(z)$.
Hence, suppose that $\tau$ is a meromorphic function on $R$ such that $\tau(z)
\cap C_z = \emptyset$ if $z \notin E_\infty$, while if $z_0 \in
E_\infty$, then $\tau(z) = O(|z-z_0|^{m+1})$.
Then $F_1^{-1}(z, \tau(z)) = (z, \tau(z)/f_1(z))$
is a holomorphic section which avoids $D$ when $z \notin E_\infty$ and
which equals $0 \in \CC \setminus \{u_0(z)\}$ when $z \in
E_\infty$.  Thus $\sigma(z) = \tau(z)/f_1(z)$ gives a holomorphic
section in the complement of $D$, so it suffices to construct $\tau$. 

As a final simplification to $D$, we note that any symmetric polynomial in
$v_0(z)$ and $v_1(z)$ gives a holomorphic function 
on $R$. So we may define
the biholomorphic map $F_2(z,w) = (z, w - (v_0(z)+v_1(z))/2)$ and let
$B = F_2(C)$.  Then $B$ has the form $B = \{(z,w) : w^2 = h(z)\}$,
where $h(z) = ((v_0(z)-v_1(z))/2)^2=((v_0(z)+v_1(z))/2)^2-v_0(z)v_1(z)$
is a holomorphic function on $R$.   Note that if $z_0 \in
E_\infty$, then  
\begin{equation}  \label{eqn:section}
F_2(z, \tau(z)) = \left(z, - \frac{u_0(z)(z-z_0)^m + u_1(z)}{2} +
O(|z-z_0|^{m+1} \right).
\end{equation}
So it suffices to construct a meromorphic function, 
$\mu$, whose graph avoids $B$ except at points
$z_0 \in E_\infty$, at which points $\mu$ should have a development as in
the second coordinate of (\ref{eqn:section}).  

To obtain the function $\mu$, we will first define a  holomorphic 
function on $B$, which we will write in the form $g(z,\sqrt{h(z)})$. We
then define
\begin{equation} \label{eqn:mu}
\mu(z) = \sqrt{h(z)} \;\frac{g(z, \sqrt{h(z)}) + g(z, -\sqrt{h(z)})}
{g(z, \sqrt{h(z)}) - g(z, -\sqrt{h(z)})}.  
\end{equation}
Note that for each $z$, $\mu$ is independent of the choice of
$\pm\sqrt{h(z)}$, so $\mu$ is a well-defined  meromorphic function 
on $R$ (by the Big Picard's theorem for example).
We will construct $g$ to ensure that $\mu$ avoids $B$
except at points in $E_\infty$, where $\mu$ has the form dictated by
(\ref{eqn:section}).  

To this end, let $E_0 = \{z: \#D_z = 1\} = \{z: 0 \in B_z\} = \{z:
h(z) = 0\}$.  We make the following requirements for the holomorphic 
function $g$ on $B$:
\begin{itemize}
\item[(a)] If $z_0 \in E_0$, then $g(z, \sqrt{h(z)}) = 1 + \sqrt{h(z)}
+ O(\sqrt{|zh(z)|})$. \\
\item[(b)] If $z_0 \in E_\infty$, then taking $\sqrt{h(z)} =
(v_0(z)-v_1(z))/2$, we have $g(z, \sqrt{h(z)}) = 1 + z^m + O(|z-z_0|^{m+1})$ 
and $g(z, -\sqrt{h(z)}) = az^m + O(|z-z_0|^{m+1})$, where $a =
u_0(z_0)/u_1(z_0)$.  \\
\item[(c)] If $z_0 \notin (E_0 \cup E_\infty)$, then $g(z,
\pm\sqrt{h(z)}) \neq 0$.  
\end{itemize}

To obtain such a function $g$, let $\Psi:X \rightarrow B$ be the
normalization of $B$, where $X$ is a noncompact Riemann surface (see
e.g. \cite{chirka}).  By \cite[theorems 26.3, 26.5]{forster}, the
classical theorems of Weierstrass and Mittag-Leffler are valid on
$X$.  Hence given finite holomorphic jets at a discrete set of points
in $X$, there is a function holomorphic on $X$ agreeing with the given
jets and having no other zeros.  This can be done for instance by
first finding a Weierstrass function, $f$, having the prescribed zeros, 
dividing all of the jet data by the jets of $f$ at the appropriate
points, then taking log of these new jets and interpolating to get a
function $r$ having these log jets.  Then $f \exp(r)$ has the desired
jets and no other zeros.  Since $\Psi$ is locally biholomorphic
except where $h=0$, the data in parts (b) and (c) transfer directly
via $\Psi$ to data on $X$.  Hence we need focus only on part (a).

Suppose $z_0 \in E_0$.  Without loss we may take $z_0 = 0$.  Suppose
that $h(z) = z^m h_1(z)$ with $m$ odd and $h_1(0) \neq 0$.  Then after
a local change of coordinates in $X$ and a choice of $\sqrt{h_1(0)}$,
$\Psi$ has the form $\Psi(x) = (x^2, x^m \sqrt{h_1(x^2)})$; i.e.,
$x=\sqrt{z}$.  In this coordinate system, $\Psi^{-1}(z, \pm
\sqrt{h(z)}) = \pm x = \pm \sqrt{z}$, and we may take $\sqrt{h(z)} = x^m
\sqrt{h_1(x^2)}$.  Hence, if $g_1$ is holomorphic in a
neighborhood of $0$ in $X$ and satisfies 
$g_1(x) = 1 + x^m \sqrt{h_1(0)} + O(|x|^{m+1})$, then $g = g_1 \circ
\Psi^{-1}$ has the  development given in (a) in a deleted
neighborhood of 0.  

Suppose next that $h(z) = z^{2m} h_1(z)$ with $h(0) \neq 0$.  In this
case, $B$ consists of two distinct components which intersect at the
origin, and hence the origin lifts to two distinct points $x_0, x_1$
in $X$.  Making a choice of $\sqrt{h_1(0)}$, we see that after a local
change of coordinates near $x_j$, $\Psi$ has the form
$$ \Psi(x) = (x-x_j, (-1)^j (x-x_j)^m \sqrt{h_1(x-x_j)}).  $$
Hence $\Psi^{-1}(z, (-1)^j z^m \sqrt{h_1(z)}) = z + x_j$.  As before,
in $X$, if $g_1$ is holomorphic in a neighborhood of $x_j$ and
satisfies $g_1(x) = 1 + (-1)^j (x-x_j)^m \sqrt{h_1(0)} + O(|x-x_j|^{m+1})$,
then $g = g_1 \circ \Psi^{-1}$ satisfies the development in (a).  

Hence the requirement in part (a) can be satisfied by specifying a
discrete set of finite jets in $X$.  As mentioned above, there exists
a function $g_1$ holomorphic on $X$ agreeing with the given finite
jets and having no other zeros.  Then $g = g_1 \circ \Psi^{-1}$
satisfies (b) and (c) and satisfies (a) in a deleted neighborhood of
each $z_0 \in E_0$.  

\medskip
Observe that for $z$ outside $E_0 \cup E_\infty$, $\mu(z) \neq \pm \sqrt{h(z)}$
since $g(z, \pm \sqrt{h(z)}) \neq 0$ for either
choice of $\pm$ there.  Hence to finish
the proof we need show only that $\mu$ is meromorphic and nonzero at
each $z_0 \in E_0$, and that $\mu$ has the development given by
(\ref{eqn:section}) at points $z_0 \in E_\infty$.  

Suppose first that $z_0 \in E_0$.  As before, we may assume $z_0 = 0$
and that $h(z) = z^m h_1(z)$ with $h_1(0) \neq 0$ and $m \geq 1$.
Recalling that $\mu$ is independent of the choice of $\pm
\sqrt{h(z)}$, we have 
\begin{align*}
\mu(z) &= \sqrt{h} \; \frac{(1 + \sqrt{h}) + (1 - \sqrt{h}) + O(|z|^{(m+1)/2})}
{(1+\sqrt{h}) - (1-\sqrt{h}) + O(|z|^{(m+1)/2})} \\
& = 1 + O(|z|^{1/2}).
\end{align*}
Hence $\mu$ has a removable singularity at $z_0 = 0$ and is nonzero
there.  

Finally, suppose that $z_0 \in E_\infty$, and recall that in this
case, $h(z) \neq 0$, and $\sqrt{h(z)} = (v_0(z)-v_1(z))/2$, where,
assuming that $z_0 = 0$, we have $v_0(z) = z^m u_0(z)$ and $v_1(z) =
u_1(z)$, $u_j(0) \neq 0$, and the development in part (b) holds.   Then 
\begin{align*} 
\mu(z) &= \left(\frac{v_0-v_1}{2}\right) 
\frac{1+(1+a)z^m + O(|z|^{m+1})}{1+(1-a)z^m + O(|z|^{m+1})} \\
&=\left(\frac{u_0 z^m - u_1}{2}\right) (1 + 2a z^m + O(|z|^{m+1})) \\
&=-\frac{u_0 z^m + u_1}{2} + u_0 z^m - a u_1 z^m + O(|z|^{m+1}).
\end{align*}
Since $a = u_0(z_0)/u_1(z_0)$, we see that this last expression has
the same form as the second coordinate in (\ref{eqn:section}).  So
$\mu$ has all the desired properties, and
hence applying $F_2^{-1}$ and $F_1^{-1}$ gives a holomorphic section
avoiding $D$, as claimed.  $\BOX$

\subsection{Normal form and dominating map}

In this section we prove corollary~\ref{cor:normal} and use this to
prove theorem~\ref{thm:dominate}.   The proof of
corollary~\ref{cor:normal} is essentially a matter of changing coordinates
in $R \times \PP^1$ so that the section provided by
theorem~\ref{thm:section} becomes the infinity section in the new
coordinates, then averaging to put the resulting section in the
desired form.  We provide the proof for completeness.  

\bigskip

{\bf Proof of corollary~\ref{cor:normal}:}  Let $D$ be a double
section over $R$.  We first apply a map which is biholomorphic on $R
\times \PP^1$ and which maps $D$ into $R \times \CC$.  To construct
such a map, note that by theorem~\ref{thm:section}, there exists a 
function $\sigma$ meromorphic on $R$ such that the graph, $\Sigma$, of
$\sigma$ in $R \times \PP^1$ is disjoint from $D$.  If $\sigma$ is
holomorphic, we apply the map $M(z,w) = (z, 1/(w-\sigma(z)))$, which is
biholomorphic on $R \times \PP^1$ and which takes $\Sigma$ to the
infinity section and hence satisfies $M(D) \subset R \times \CC$.  

Otherwise, $\sigma$ is meromorphic but not holomorphic.  Choose a
function $\psi$ holomorphic on $R$ such that $\psi$ is zero precisely
when $\sigma$ has either a zero or a pole, with the zeros of $\psi$
chosen so that $1/\psi$ has the same principal part as $1/\sigma$ if
$\sigma$ has a 0 and $1/\psi$ has the same principal part as $\sigma$
if $\sigma$ has a pole.  Let $\phi(z) = (-1/\sigma(z)) + (1/\psi(z))$,
and let
$$ M(z,w) = \frac{\phi(z) w + 1}{w - \sigma(z)}. $$
If $\sigma(z_0) \in \CC$, then $\phi(z_0) \in \CC$ with $\phi(z_0) \neq
-1/\sigma(z_0)$, and hence $M(z_0, \cdot)$ is a nondegenerate Mobius
transformation taking $\sigma(z_0)$ to $\infty$. If $\sigma$ has a pole
of order $m$ at $z_0$, then $\phi$ has a pole of the same order at $z_0$,
so multiplying numerator and denominator by $(z-z_0)^m$, we see that
$M(z_0, w) = cw$ for some $c \neq 0$.  Hence $M$ is biholomorphic and
$M(\Sigma)$ is the infinity section, so $M(D)$ is contained in $R
\times \CC$.  

Now, near a generic point in $R$, $M(D)$ is the union of the graphs of
two holomorphic functions $u_0$ and $u_1$, and the symmetric function
$u_0 + u_1$ is holomorphic on $R$.  Hence the map
$$ N(z,w) = \left(z, w - \frac{u_0(z)+u_1(z)}{2}\right)$$
is biholomorphic on $R \times \PP^1$, and taking $h(z) =
(u_0(z)-u_1(z))^2/4$, we have $N(M(D)) = \{(z,w): w^2 = h(z)\}$.  $\BOX$

\bigskip

Finally, we are ready for the construction of a dominating map into
the complement of a double section.

\bigskip
{\bf Proof of theorem~\ref{thm:dominate}:}  Let $D$ be a double
section in $R \times \PP^1$.  By the corollary, we may assume that $D
= \{(z,w): w^2 = h(z)\}$ for some $h$ holomorphic on $R$.  Let 
$$ H(z,w) = \sqrt{h(z)} \; \frac{\exp(\sqrt{h(z)}w) + \exp(-\sqrt{h(z)}w)}
{\exp(\sqrt{h(z)}w) - \exp(-\sqrt{h(z)}w)}.$$
Note that $H$ is independent of the choice of $\pm \sqrt{h(z)}$.
Moreover, $H$ is simply the Mobius transformation $x \mapsto
(x+1)/(x-1)$ applied to $\exp(2\sqrt{h(z)}w)$, then multiplied by
$\sqrt{h(z)}$.  In particular, if $h(z_0) \neq 0$, then $H(z_0,
\cdot)$ maps $\CC$ into $\PP^1 \setminus \{\pm \sqrt{h(z)}\}$, and $F$
is holomorphic except possibly where $h(z) = 0$.  

Suppose $h(z_0) = 0$, and without loss take $z_0 = 0$.  Then $h(z) =
z^m h_1(z)$ for some $m \geq 1$, $h_1(0) \neq 0$.   Near $0$ we have
\begin{align*}
H(z,w) &= \sqrt{h(z)} \; 
\frac{(1+\sqrt{h(z)}w) + (1-\sqrt{h(z)}w) + O(|z|^{(m+1)/2})}
     {(1+\sqrt{h(z)}w) - (1-\sqrt{h(z)}w) + O(|z|^{(m+1)/2})} \\
& = \frac{1 + O(|z|^{1/2})}{w + O(|z|^{(m+1)/2})}.
\end{align*}
Hence for each fixed $w\neq 0$, $H(\cdot, w)$ has a removable singularity at
$0$, and $H(0, w) = 1/w$.  So $H$ is holomorphic in each variable
separately for $w\neq 0$, hence holomorphic by Hartogs' theorem
as we have seen that $H$ is holomorphic
for $z\neq 0$.  Moreover, for each
fixed $z$, $H$ has the form of an exponential map followed by a Mobius
transformation if $D$ has two points in the fiber above $z$, and $H$
has the form of a Mobius transformation if $D$ has one point in the
fiber above $z$.  $\BOX$

\subsection{Jet prescriptions}

In this section we prove theorem~\ref{thm:jet}.  So, let $D$ be a
double section over $R$ and let $E$ and $P_c$ be as in the statement
of the theorem.  As before, we can interpolate on $R$, so there exists
a function $\Psi$ meromorphic on $R$ such that if $c \in E$, then the
Laurent expansion of $\Psi$ at $c$ agrees with $P_c$ to order
$\deg(P_c)$.  We may assume also that the graph of $\Psi$ is not
contained in a component of $D$.  As in the proof of
corollary~\ref{cor:normal}, we can make a global biholomorphic change
of variables so that the graph of $\Psi$ becomes the infinity section
in the new coordinates.  Then the requirement of agreeing with $P_c$
becomes the requirement of having a pole of at least some given order
$n_c$ at each $c \in E$.   

Let $p$ be a holomorphic function on $R$ 
with a zero of order $n_c$ at each $c \in E$ and no
other zeros, and let $\Phi(z,w) = (z, p(z)w)$.  Then $\Phi(D)$ is a
double section over $R$.  We may assume that $n_c$ is chosen large
enough that if $c \in E$, then $\Phi(D)$ intersects the fiber over $c$
in the single point 0.  By theorem~\ref{thm:dominate}, there exists
a dominating map, $G$, from $R\times \CC$ onto $(R \times \PP^1) \setminus
\Phi(D)$ of the form $G(z,w) = (z, L(z,w))$ as described in that
theorem.  

To finish the proof of the current theorem, let $F(z,w) =
\Phi^{-1}G(z,w) = (z, L(z,w)/p(z))$.  If $c \notin E$, then $F(c,
\cdot)$ is easily seen to map $\CC$ onto $\PP^1\setminus D_c$,
while if $c \in E$ and $\sigma$ is holomorphic in a neighborhood of
$c$, then $L(z, \sigma(z))$ is meromorphic and nonzero in a
neighborhood of $c$, hence $F(z, \sigma(z))$ has a pole of order at
least $n_c$ at $c$.  Thus $F$ satisfies all the requirements. $\BOX$

\section{The Jacobian fibration route}

In this section we utilize and extend the theory of Jacobian
fibrations introduced by Kodaira in \cite{Kod} to obtain the simultaneous
uniformization result of theorem~\ref{thm:Jac}.
Elementary lemmas and arguments will be given
in an attempt at wider accessibility. Nevertheless,
details on many of the background material can be found in \cite{BPV},
especially section V.9, and the references therein. The second author
would like to dedicate this section to the memory of K. Kodaira, whom
he came close to meeting.

Our starting setup is a $\PP^1$-bundle $p:P\rightarrow R$ with a 
distinguished double section $D$. A fiber 
here is just $P_z=\PP^1$ marked with
an effective divisor $D_z$ of degree $2$. In
the spirit of Kodaira,
there are two simple ways to attach a nontrivial Jacobian, or
more precisely a Picard variety, to this marked fiber by
turning it into a singular curve. One way is to take
the Jacobian of the singular curve 
obtained by ``collapsing, to first order,'' the
divisor $D_z$ into one point, resulting in a reduced and irreducible 
curve isomorphic to either
a nodal or a cuspidal cubic curve in $\PP^2$. This can be done via
the linear system generated by $D_z+q_1$, $D_z+q_2$ and $3q_1$ where
$q_1$ and $q_2$ are points in $P_z$
chosen in general position with respect to $D_z$.
The other way is to add to $P_z$, which is a smooth rational
curve, another  smooth rational curve
intersecting $P_z$ at $D_z$ counting multiplicity. The combined object
can again be realized as one of two types of
plane cubic curves, but which are reducible, corresponding to the 
reduced total transform of the two types of irreducible 
plane cubic curves above blown up at their singular point. 
Although these methods are equivalent, and both will be introduced below, 
we will focus on the second one, which is easier to introduce.

\medskip

We first recall that the Picard variety of a one dimensional
projective scheme $Z$ is the space $\mbox{Pic}^0(Z)$ of analytic
line bundles with vanishing first Chern class on $Z$, 
where the first Chern class $c_1$ is defined by the connecting 
homomorphism $\delta$ in the long exact cohomology sequence over $Z$ of
the exponential sequence (i.e., the short exact sequence 
$0 \rightarrow \ZZ \rightarrow {\cal O}\rightarrow {\cal O} ^ {\ast} 
\rightarrow 0$):
$$
\begin{array}{clclclclc}
0  &\rightarrow & H ^{1 } (Z, \ZZ ) & \stackrel{i}{\rightarrow} & 
H ^{1 } ( {\cal O} _{Z}
) & \rightarrow & H ^ {1 } ( {\cal O} _{Z} ^ {\ast} ) &
\stackrel{\delta}{\rightarrow} & H ^ {2} (Z, \ZZ ) \rightarrow 0 \\
   &\ & || & & || & & || & & \\
   &\ &\mbox{a f.g. }\ZZ\mbox{-module} & & \CC^{\, g} 
& & \mbox{Pic} (Z) & & 
\end{array}
$$
\noindent
Here, $g=p_a(Z)$ is the arithmetic genus of $Z$ 
and Pic$(Z)$ the space of holomorphic line bundles
over $Z$. If $Z$ is reduced on the smooth part $Z_i^\circ$
of an irreducible component $Z_i$ of $Z_{red}$, then
every line bundle $L$ on $Z$ can be written as
${\cal O}(E+E')$ for a divisor $E=\sum n_j s_j$ on $Z_i^\circ$
($n_j\in \ZZ, s_j\in Z_i^\circ$) and a divisor $E'$ on $Z\setminus Z_i$.
If further $L$ is trivial on 
the Zariski closure of $Z\setminus Z_i$, then we can take $\delta(L)$
to be $\deg E=\sum n_j$. We note that Pic$^0(Z)$ has a natural
structure of a complex Lie group given by the tensor product of line bundles.
Also, in the case when $g=1$, $Z_i^\circ$ can be shown to be
biholomorphic to Pic$^0(Z)$ via Serre duality and the Riemann-Roch
theorem for singular curves (c.f. section V.9 of \cite{BPV}). This fact is
well understood classically for reduced plane cubic curves and 
the construction to see this is standard. However, for the benefit
of readers who are not familiar with the projective
geometry of the plane cubic curves,
an exposition of this simple but beautiful explicit classical construction 
will be given in the proof of the lemma to follow 
(see, for instance, examples 6.10.1 and 6.11.4] of \cite{Ha}). 

As a preamble, we note that all reduced plane cubic curves occur in 
Kodaira's table of singular fibers; see table 3 in Chapter V of \cite{BPV}. 
Indeed, the classical formula for the arithmetic genus $g=(d-1)(d-2)/2$, where
$d$ is the degree of the plane curve, shows that they all have $g=1$.
They are classified as Kodaira's
I$_0$, corresponding to the nonsingular elliptic curves, I$_1$ and II,
corresponding to the plane cubic curves with a node and the ones with a cusp
(i.e., rational curves with a node and ones with a cusp), I$_2$ and III,
corresponding to the reduced total transform of the 
blowup of the above two types (respectively) at their 
singular points in the plane, and IV, corresponding
to the reduced total transform of the blowup of type III at the singular point.
We remark that each type of fiber occuring in Kodaira's list just mentioned can
be regarded as a plane cubic via Riemann Roch and we assume this in the 
following proof. Of the singular plane cubic curves, only 
those of type I$_1$ and II are irreducible and they will serve
as the starting point for the proof below.

\begin{lemma}  \label{lemma:Inversion} 
Let $Z$ be a reduced plane cubic curve with
a base point $\sigma\in Z_i^\circ$, a component of the smooth
locus of $Z$. Then
$Z_i^\circ$ is biholomorphic to ${\rm Pic}^0(Z)$ via the map
$$
x \in Z_i^\circ \stackrel{f}{\mapsto} {\cal O} (x - \sigma ) 
\in {\rm Pic} ^ {0} (Z).
$$
\end{lemma}
\noindent
{\bf Proof:} We first show that $f$ is one-to-one for the case $Z$ is
an irreducible cubic curve 
with a node (type I$_1$) and the case with a cusp (type II).
Assume not, so that ${\cal O} (x - \sigma ) = {\cal O} ( x ^ {\prime} -
\sigma)$ where $x \neq x ^{\prime}$. Then 
${\cal O} (x- x ^{\prime} )$ corresponds to the
trivial line bundle over $Z$ and so $Z$ has a rational 
function with a simple pole at $x ^ {\prime}$ and a simple zero at $x$.
This gives a birational morphism from $Z$ to $\PP^1$.
But any such morphism composed with the normalization morphism for $Z$
would give an automorphism of $\PP^1$ that factors through $Z$ which
is impossible since then two different points (infinitely near in the 
cuspidal case) of $\PP^1$ would map to the same point by this automorphism.

Since the other type of reduced singular cubics
can be regarded as the reduced total transforms of the above 
two cases via the blowing up process, the same argument as above
works for them because the birational morphism
from $Z_i$ to $\PP^1$ obtained as above
would factor through the corresponding irreducible cubic curves
considered above. The case of the smooth cubic also follows the same argument
via the birational invariance of the geometric genus for smooth curves, 
which is $0$ for $\PP^1$. Hence $f$ is one-to-one for all 
the cases of plane cubics.

We next show that $f$ is onto. Again, we do this
first for the two cases of $Z$
corresponding to that of an irreducible plane cubic
with a node and that with a cusp. In these cases, we may write $Z_i^\circ$
simply as $Z^\circ$. Consider a tangent line
to an inflection point $o$ on $Z^\circ$ 
(there are exactly $9$ of these on any irreducible plane cubic, all ordinary), 
which would meet $Z$ in the divisor $3o$. Any other
divisor $H$ obtained by linear intersection would be linearly
equivalent to $3o$, and we denote this by $H\sim 3o$. Let $L\in
\mbox{Pic}^0(Z)$.  
We may write $L\otimes {\cal O}(\sigma -o)$
in the form ${\cal O}(E)$ for a divisor $E=\sum n_j s_j$ with
$n_j\in \ZZ$ and $s_j\in Z^\circ$ for all $j$ and $\sum n_j=0$.
It follows that $E=\sum n_j (s_j-o)$ also. Now, given a point $s\in Z^\circ$, 
consider the line intersecting $Z$ at the divisor $o+s+t$, where $t\in Z^\circ$
by counting multiplicity of the intersection. Then $o+s+t\sim 3o$, giving
$s-o\sim -(t-o)$. This allows us to write $E=\sum n_j(s_j-o)$ with all $n_j>0$.
We now show by induction on $\sum n_j$ that $E\sim s-o$ 
for some $s\in Z^\circ$. So suppose $\sum n_j\geq 2$ and let $s,s'$ be two of
the points $s_j$ (maybe the same) occurring in the expression for $E$. We
now consider the lines intersecting $Z$ at $s+s'+t$ and at $o+t+z$. As
before, we see $t,z\in Z^\circ$ and $(s-o)+(s'-o)\sim (z-o)$ which allows
us to reduce $\sum n_j$ by $1$ and consequently $E\sim  s-o$ 
for some $s\in Z^\circ$ via induction. 
Thus $L={\cal O}(E-\sigma +o)={\cal O}(s-\sigma)$ 
for some $s\in Z^\circ$ and so $f$ is onto in these two cases of $Z$.
The other cases of singular $Z$ now follow directly as before since
they are constructed from the above cases via blowing up. Finally,
the case of the smooth cubic also follows verbatim the argument above.
\smallskip

We have in fact exhibited the Lie group structure imposed by Pic$^0$
on the smooth part of
the cubic curves above via linear intersection.
Note that the third point in a linear intersection of a cubic above depends
holomorphically on the other two. It is now a formality to verify that $f$
is holomorphic. This concludes our proof. $\BOX$
\smallskip

We should point out that a common but more abstract and unified 
approach to show
the surjectivity of $f$, without reference to which
cubic,  is to apply the Riemann Roch theorem
and Serre duality for an embedded curve to the line sheaf $L(\sigma)$ 
where the dualizing sheaf of our curve $Z$ is just
${\cal O}_Z$ by the adjunction formula (see sections II.1-6 of \cite{BPV}).

\bigskip

We will introduce the Jacobian fibration for our 
$\PP^1$-bundle $p:P\rightarrow R$ endowed
with the double section $D$
as the Jacobian fibration Jac$(\bar p)$ of an associated family 
$\bar p: \bar P\rightarrow R$ of singular plane cubic curves over $R$. 
The most direct way to do this is to realize $p$ as a family of conics
(curves of degree two) in a $\PP^2$ bundle over $R$ via the canonical
embedding given by $D$ and attach
to it the family of lines specified by $D$. More specifically:

\begin{prop}
With the setup as above, ${\cal E}=p_*{\cal O}(D)$ is a locally free sheaf of
rank $3$ over $R$. The natural morphism 
$\Phi_{D/R}:P\rightarrow \PP({\cal E})$ is an embedding with
${\cal O}(D)=(\Phi_{D/R})^*{\cal O}_{\PP({\cal E})}(1)$.
Furthermore, ${\cal O}_{\PP({\cal E})}(1)$ has a section $\nu$ 
such that $(\Phi_{D/R})^*(\nu)=D$.
\end{prop}

\noindent
{\bf Proof:} Since ${\cal O}(D)_z$ is a divisor of degree $2$ on 
$P_z=\PP^1$, it is simply ${\cal O}_{\PP^1}(2)$ for all $z$. So
$h^0({\cal O}(D)_z):=\dim H^0({\cal O}(D)_z)=3$ is independent of $z$
and our first claim follows from Grauert's theorem, theorem 12.9
of \cite[\S III]{Ha}, as $p$ is flat. 
Now, the base change property of the $0$-th direct
image map from the same theorem coupled with Nakayama's lemma shows
that the natural map $p^*{\cal E}\rightarrow {\cal O}(D)$ is surjective
as ${\cal O}(D)_z$ is generated by global sections. Hence the dual of
this natural map projectivizes to a morphism 
$\Phi_{D/R}:P\rightarrow \PP({\cal E})$ with 
${\cal O}(D)=(\Phi_{D/R})^*{\cal O}_{\PP({\cal E})}(1)$.
This gives an embedding on each fiber as
${\cal O}(D)_z={\cal O}_{\PP^1}(2)$ is very ample. 
Hence $\Phi_{D/R}$ is an embedding and by definition, $D$ corresponds
to a section of ${\cal O}_{\PP({\cal E})}(1)$. $\BOX$

\medskip

We may now take $\bar P$ to be the union of $\Phi_{D/R}(P)$ and 
the two dimensional surface corresponding to $(\nu)$ in ${\PP({\cal E})}$
and take $\bar p$ to be the bundle projection of $ \PP({\cal E})$ restricted
to $\bar P$.  By construction, $\bar p$ is easily verified to give a flat 
family of reduced schemes. (A less direct way to see this is to show that
there are no embedded points in codimension two
as $\bar p$ is flat from the general theory,
such as proposition~9.7 of \cite[\S III]{Ha}. But this is clear from the fact
that all reduced plane cubic curves have arithmetic genus $g=1$, and an
embedded
point would increase the arithmetic genus, which would then contradict the
fact that arithmetic genus stays constant in a flat family.)
Note that each fiber of $\bar p$ as a scheme is a singular cubic curve with
two irreducible components consisting of a line and a conic 
intersecting precisely at $\Phi_{D/R}(D)$.
Note also that $p=\bar p\circ \Phi_{D/R}$ by construction.

The exponential sequence on $\bar P$ leads to the exact sequence
$$
0\rightarrow \bar p_*\ZZ \rightarrow \bar p_*{\cal O}_{\bar P}
\rightarrow \bar p_*{\cal O}^*_{\bar P} \rightarrow 
\bar p_{*1}\ZZ \rightarrow \bar p_{*1}{\cal O}_{\bar P}
$$
where the zeroth and the first direct image of ${\cal O}_{\bar P}$
are locally free sheaves of rank one by Gauert's theorem as $\bar p$
is flat and $g$ is constant. In fact, the natural structure map 
${\cal O}_R\rightarrow \bar p_*{\cal O}_{\bar P}$ generates the stalk
everywhere by the base change property so that 
$\bar p_*{\cal O}_{\bar P}={\cal O}_R$ forcing
$\bar p_*{\cal O}^*_{\bar P}={\cal O}^*_R$ and $\bar p_*\ZZ=\ZZ_R$ as
both are subsheaves. It follows that the beginning of the above sequence
is nothing but the exponential sequence over $R$ so that we can regard
$\bar p _{\ast 1} \ZZ$ as a subsheaf of $\bar p_{\ast 1}{\cal O}_{\bar P}$. 
The latter, being locally free, 
is the sheaf of sections of a line bundle $L$ over
$R$ and so includes $\bar p _{\ast 1} \ZZ$ as a subset.

To construct our
Jacobian fibration, we first form the  ${\cal O} _{R}$ module 
$${\cal J}\! ac ( \bar p ) = {\bar p}_{\ast 1} {\cal O} _{\bar P}  /
\bar p _{\ast 1} \ZZ$$ over $R$. We would like to construct a holomorphic
quotient space of $L$ by $ \bar p _{\ast 1} \ZZ$ so that
${\cal J} \! ac ( \bar p )$ corresponds
to the sheaf of sections of 
	$$\mbox{Jac}( \bar p ):= L / \bar p _{\ast 1} \ZZ,$$ 
which would then be a holomorphic fibration of complex Lie
groups with a zero section and which we will define as our Jacobian fibration.
This quotient space is certainly possible in a neighborhood of any point
where the double section $D$ consists of two distinct points since
$\bar p _{\ast 1} \ZZ$ is a closed submanifold of $L$ generated over $\ZZ$
by a nonvanishing section of $L$ there.
To see that $\bar p _{\ast 1} \ZZ$ is actually closed everywhere in $L$
it is necessary only to observe that in a neighborhood of 
any point where $D$ becomes a double point, 
the double valued generator over $\ZZ$ of 
$\bar p _{\ast 1} \ZZ$ becomes a nonvanishing
meromorphic single valued generator after pulling back
the fibration $p$ and hence $\bar p$ via the normalization of $D$, 
showing that the original double valued generator approaches $\infty$
at this double point. Hence it is
straightforward to construct the quotient 
space (see \cite{Kod}, compare also \cite[V.9]{BPV} and 
Grothendieck \cite{Gro}, who showed that Jac exists as
a scheme for any flat family) where the quotient map is a 
local isomorphism.
This quotient is thus a complex manifold
with a canonical section over $R$ coming from the 
zero section of $L$. In fact, we have the following 
proposition, which should be compared
with proposition 9.1 of \cite[V.9]{BPV}.

\begin{prop}   \label{prop:L3}
With the setup as above, suppose that $p$ has a holomorphic
section $\sigma : \ R \rightarrow P$ disjoint from $D$.
Then there is a canonical fiber-preserving
isomorphism $h$ from {\rm Jac}$( \bar p )$ onto $P^{\circ}=P\setminus D$ 
mapping the zero-section in {\rm Jac}$( \bar p )$ onto $\sigma (R)$.
\end{prop}

\noindent
{\bf Proof:} The base change property shows that 
Pic$^0(\bar P_z)={\rm Jac}( \bar p )_z$ thanks to the flatness of
$\bar p$ and the constancy of the arithmetic genus of its fibers.
Let $P'$ be the irreducible component of $\bar P$ corresponding to the
reduced divisor $(\nu)$ in our construction of $\bar P$ above. Then
$P^\circ$ is isomorphic to $\bar P\setminus P'$ under $\Phi_{D/R}$
where $P^\circ_z$ is identified with
the component of the smooth part of $\bar P_z$ containing 
$\sigma(z)$. Consider the map $J:P^\circ\rightarrow {\rm Jac}( \bar p )$
given on each fiber by 
$$x\in P^\circ_z\mapsto {\cal O}(x-\sigma(z))\in {\rm Pic}^0(\bar P_z).$$
Lemma~\ref{lemma:Inversion} 
shows that this is a bijective map that is holomorphic
in the fiber direction. It is trivial to verify that $J$ maps holomorphic
sections to the same locally. Hence $J$ is biholomorphic 
and we may set $h=J^{-1}$. $\BOX$
\bigskip

It is important to note that sections of ${\rm Jac}( \bar p )$, and hence of
${\cal O}(L)= \bar p_{\ast 1}{\cal O}_{\bar P}$, act via addition on
$P\setminus D$ as soon as 
$P\setminus D$ has a section by this proposition.

Theorem~\ref{thm:Jac} is a corollary of this proposition.  Moreover,
using this proposition, we now give a second proof of
theorem~\ref{thm:section} following the  proof of
Lemma~{3.8} in \cite{bl} with some needed clarification.
Theorem~\ref{thm:dominate} then clearly follows  by
a final application of the proposition. 

\bigskip
\noindent
{\bf Second proof of theorem~\ref{thm:section}:} 
With the setup as above, the hypotheses of theorem~\ref{thm:section}
allow us to assume that $R$ is noncompact and $P=R\times \PP^1$ with 
$p$ being the projection to the first factor. 
We will exploit the fact that $R$, being Stein, contains 
an one-dimensional simplicial
complex $K$ which is a strong deformation
retract of $R$ (see \cite[\S V.4.5]{GR}).
\smallskip

	We first assume that $D$ has no double point at any
fibre, which is equivalent to assuming that 
$\Gamma=\bar p_{\ast 1}{\cal O}_{\bar P}$ is a locally
constant sheaf, or equivalently, a local system. 
Then the pair $(K, \Gamma|_K)$ is a deformation 
retract of $(R,\Gamma)$. Hence $H^2(R, \Gamma)=H^2(K,\Gamma)$, which
is zero as it is isomorphic to the simplicial sheaf cohomology of a
simplicial complex (since $\Gamma$ restricted to any simplex is a
constant sheaf and hence has trivial cohomology, leading to the
degeneration at $E_1$ of the double complex comparing simplicial and
ordinary sheaf cohomology) which has dimension less than two. Since
the cohomology groups of $L$ are trivial as $R$ is Stein, the long
exact sequence of the short exact sequence
\begin{equation}\label{J}
0\ra \Gamma \ra L \ra {\rm Jac}( \bar p ) \ra 0
\end{equation}
(as sheaves) shows that $H^1(R,{\rm Jac}( \bar p ))=H^2(R,\Gamma)$.
Hence $H^1(R,{\rm Jac}( \bar p ))=0$.

	Now, every point on
${R}$ admits a neighborhood with a section of $p$ avoiding $D$.
We choose a locally finite good covering of $R$ by distinct open sets $U _{1} ,
U _{2} ,...,$ with sections $\tau _{1} , \tau _{2} ,...$ of 
$p$ over ${U _{1}}, {U _{2}} ,...,$ respectively.
By proposition~\ref{prop:L3}, there is a unique section 
$\tau _{ij} ^ {\prime} \in H ^{0} (U _{ij}, {{\rm Jac}( \bar p )} )$ such that
$\tau _{i} + \tau _{ij} ^ {\prime} = \tau _{j}$ on $U _{i} \cap U _{j}$
for all $i,j$. It is clear that 
$\{\tau_{ij}^{\prime}\}$ satisfies the cocycle condition. Hence, so does
$\{-\tau_{ij}^{\prime}\}$, which therefore defines an element in 
$H ^{1 } ( \{ U _{i} \} , {\rm Jac}( \bar p ))=
H ^{1 } ( R , {\rm Jac}( \bar p ))$ (where the last equality follows
from Leray's theorem \cite[\S B.3.3]{GR}). 
Hence $\{-\tau_{ij}^{\prime}\}$ is a coboundary and
one can find holomorphic sections $\tau _{i} ^ {\prime} 
\in H ^{0} (U _{i} , {{\rm Jac}( \bar p )} )$ 
such that $\tau _{i} ^ {\prime } -
\tau _{ij} ^ {\prime} = \tau _{j} ^ {\prime}$.  Then $\tau _{i} + \tau _{i} ^
{\prime} = \tau _{j} + \tau _{j} ^ {\prime}$ on $U_i\cap U_j$ for all $i,j$.
This gives rise to a global section of $p$ avoiding $D$. 
\smallskip

	In the general case, the problem again reduces to the vanishing
of $H^2(R,\Gamma)$, which allows us to lift the $\tau_{ij}'s$ as above
to sections of $L$ (by assuming without loss that
$D$ consists of two disjoint components over
$U_i\cap U_j$) that satisfy the cocycle condition and thus giving
an element of $H^1(\{U_i\}, L)$, which is zero as $R$ is Stein. This 
can be done via an elementary combinatorial argument using $K$ (see
\cite{Lu}), but we offer here an alternative to $H^2(R,\Gamma)$:

let $E_0$ be the set of points above which
$D$ is a double point. Then $E_0$ is a closed discrete subset of $R$.
Note that the closure $A$ in $R$ of a 
component of $R\setminus K$, being contractible modulo its boundary
$\partial A$ in $R$, is homeomorphic to either the complement of
$\DD$ with some identifications on the boundary circle
or to the upper half plane in $\CC$. It follows that any closed discrete 
subset of $A$ is contained in $\partial A$
union a simple (without self-intersection) path in $A\setminus \partial A$
connecting $\partial A$ with the complement of every compact subset of $A$
(by using, for example, the discrete set of circles in $\CC$
centered at the origin which together contain the discrete 
subset in $\CC$ corresponding to $E_0\cap A$). Hence we may define 
a one-dimensional simplicial complex $J$ containing $K$ and $E_0$ as
the union of $K$ with a set of such paths, one path for each
component of $R\setminus K$. Then $\Gamma$ is
locally constant on $R':=R\setminus J$ and each component of $U\setminus J$
is contractible for every deformation retract neighborhood $U$ of
$J$ in $R$. We take such a neighborhood $U$ small enough so that it has 
a covering by distinct open subsets $U_1,U_2, \dots$ of $U$
without triple
intersections and with sections $t_1, t_2, \dots$ of $p$ over
$U_1,U_2, \dots$, respectively, all avoiding $D$.
We may further stipulate that $D$ consists of two disjoint components over
each double intersection $U_{ij},\ i\neq j.$ Then $\Gamma$ is locally
constant on $U_{ij}$, hence constant, giving $H^1(U_{ij},\Gamma)=0$.
It follows that the section 
$t_{ij}^{\prime}$ of ${\rm Jac}( \bar p )$ over $U_{ij}$
constructed via proposition~\ref{prop:L3} to satisfy 
$t_i+t_{ij}^{\prime}=t_j$ lifts to a section $\hat t_{ij}$ of $L$ over $U_{ij}$
by the long exact sequence of equation~\ref{J}. 
As $\{\hat t_{ij}\}$ trivially satisfy the 
cocycle condition and as $H^1(U, L)=0$, there exist sections 
$\hat t_i\in H^0(U_i,L)$ such that $\hat t_i-\hat t_{ij}=\hat t_j$.
Then $\hat t_i+t_i=\hat t_j+t_j$ on $U_{ij}$ for all $i\neq j$. This
gives rise to a global section of $p$ avoiding $D$ over $U$.

	Finally, $R$ is covered by the open sets $R'$ and $U$, both of 
which admits sections of $p$ as shown above, and $R'\cap U$
consists of contractible sets by construction,  Hence the same argument
as in the last paragraph gives a global section of $p$ 
avoiding $D$ as required.
$\BOX$

\medskip

We note that this proof is very similar to the approach
of the Oka-Grauert principle for principal (and associated) bundles 
over Stein spaces (see \cite{GR} for example).  However, the current
proof applies to a
fibration of Lie groups which is not necessarily smooth.

\section{Fibered surfaces}

In this section, we will give conditions under which a nonhyperbolically 
fibered complex surface is dominable by $\CC^2$ by proving 
theorem~\ref{thm:surf}.

	Recall that a nonhyperbolically fibered surface
$X$ is given by a holomorphic map $f$ from $X$ onto a complex curve 
$C$ whose generic fiber is nonhyperbolic such that $f$
extends to a proper map from a surface $Y$ to $C$,
where $Y$ contains $X$ as a Zariski open subset. 
We will also use $f$ to denote the extended
map by abuse. 

	Recall also that each fiber $X_s$ is
an effective divisor on $X$ as the pull back of the divisor $s\in C$.
We write $X _{s} = \sum n _{i} C _{i}$ 
where $C _{i}$ is the $i$-th irreducible component of 
$(X_{s})_{red}$ and $n _{i} -1$ is the 
generic vanishing order of $d f$ on $C_{i}$.
The positive integer coefficient $n_i$ is called the multiplicity
of the $i$-th component.
The multiplicity of a fiber $X _{s} = \sum n _{i} C _{i}$ is defined
as the greatest common divisor $n_s$ of $\{ n _{i} \}$.  A fiber 
$X_{s}$ with $n_s>1$ is called a multiple fiber.

\medskip

	Vital to our problem is the minimum multiplicity in the fiber $X_s$
defined by $m_s=\min\{ n _{i} \}$. As we will soon see,
this gives the correct orbifold
structure on $C$ via the $\QQ$-divisor $A=\sum_s(1-1/m_s)s$ in $C$
for hyperbolicity questions.
Note that $X_s$ is a Zariski open subset of $Y_s$ and we will assume
for clarity that the map $f$ in a neighborhood of each point $y\in Y_s$
is given by an equation of the form
\begin{equation}\label{eqn:cd}
x_1^{n_{i_1}} x_2^{n_{i_2}}=t
\end{equation}
for a local coordinate system $(x_1, x_2)$ centered at $y$ and a
local coordinate $t$ centered at $s$. We can always do this by
further blowing up $Y$ if necessary (see theorem V.3.9 in \cite{Ha}).
We will first prove that (c) implies (b) in theorem~\ref{thm:surf}.

\subsection{Ahlfors' lemma for generalized orbifold maps}
\begin{definition} Let $(C,A)$ be a one dimensional complex orbifold, where
the $\QQ$-divisor $A=\sum_s(1-1/m_s)s$ on the complex curve $C$ specifies 
the orbifold structure as above and where the positive integer $m_s$ is 
identically $1$ except on a discrete subset of $C$. An orbifold metric
on $(C,A)$ is given by a $(1,1)$-form $\rho$ on 
$C\setminus A_{red}$ which, in a coordinate neighborhood $U(s)$ 
of each $s\in C$, has the form $$|t|^{-2(1-\frac{1}{m_s})}\mu$$
for a  smooth positive $(1,1)$-form $\mu$ on $U$ with the coordinate $t$ 
centered at $s$. 

\smallskip

A holomorphic map $f$ from a complex manifold $X$ to
$C$ is called a Z-orbifold map to $(C,A)$ if $f^*\rho$ 
extends to a continuous $(1,1)$-form over its singularities (if it has any)
on $X$ for some (and hence any) orbifold metric $\rho$ on 
$(C,A)$.
\end{definition}

A example of a Z-orbifold map is given by the fibering map $f$ of a fibered
surface $X$ over a curve $C$ with the minimum multiplicity in the fibers
defining the orbifold structure on $C$. One can deduce this directly from 
the local coordinate picture given by equation~\ref{eqn:cd}.
Composing $f$ with any holomorphic
map gives of course again a Z-orbifold map. If $(C,A)$ is hyperbolic,
meaning that it is an orbifold quotient of the unit disk $\DD$ by automorphisms
of $\DD$, then it inherits a canonical orbifold metric from the Poincar\'e
metric $\hat\rho$ on $\DD$, which is invariant under automorphisms of $\DD$.
We will call this orbifold metric the 
orbifold Poincar\'e metric and denote it by $\rho$.
It is determined as an orbifold metric on $(C,A)$ by the equation 
$$dd^c\log \rho=\rho$$
just as the Poincar\'e metric on $\DD$. With this in mind, the 
proof of the following version of Ahlfors' lemma follows verbatim that
of the usual one and is left to the reader.

\begin{theorem} \label{Ahlfors}
Denote the Poincar\'e metric on $\DD$ by $\hat\rho$ and the orbifold
Poincar\'e metric on
a hyperbolic one-dimensional complex orbifold $(C,A)$ by $\rho$. Then any 
Z-orbifold map $h$ from $\DD$ to $(C,A)$ satisfies the distance 
decreasing property in the form 
$$h^*\rho\leq \hat\rho$$
with equality at any point in $C\setminus A$ characterizing $h$ to be
the orbifold quotient map from $\DD$.
\end{theorem}

By applying this to larger and larger disks in $\CC$, we see that any 
Z-orbifold map from $\CC$ to a hyperbolic $(C,A)$ must be constant. From 
the classical uniformization theorem 
for orbifold curves, $(C,A)$ is not hyperbolic (meaning  that it
is not uniformized by the disk) precisely when $C$ is quasiprojective
and $\chi(C,A)\geq 0$ (see \cite{bl}, which explains this well known
fact from \cite[IV 9.12]{FK}). 
So if the latter condition does not hold, then
any holomorphic image of $\CC$ in a surface $X$ fibered over the 
orbifold $(C,A)$ must lie in a fiber.
This shows that condition (c) of theorem~\ref{thm:surf}
implies condition (b).

\subsection{Proof of theorem~\ref{thm:surf}:}
It is easy to deduce (c) from (a),
that if there is a dominating map $F:\CC^2 \ra X$, then there is
also a holomorphic image of $\CC$ which is Zariski dense:
First we may assume that the Jacobian of $F$ is non-zero at
the origin.  Defining $g:\CC \ra \CC^2$ by $h(z) = (\sin(2 \pi z),
\sin(2\pi z^2))$, we see that $g(n)=(0,0)$ with corresponding tangent
direction $(2\pi, 4\pi n)$ for each $n \in \ZZ$.  Taking $F \circ g$, we
obtain a holomorphic image of $\CC$ with an infinite number of tangent
directions at one point, which implies that the image is Zariski dense.
\bigskip

To complete the proof, we need to deduce (a), the existence of a 
dominating map from $\CC^2$, from condition (b), the nonhyperbolicity
of $(C,A)$. The latter implies, by the uniformization theorem, that
$(C,A)$ is the orbifold quotient of $\CC$ or $\PP^1$ unless $C=\PP^1$ with 
$\#\{A_{red}\}\leq 2$. In the ``unless'' case, by
removing up to two fibers in $X$, we may assume that $C=\CC^*$ with
$A=\emptyset$ to construct our dominating map. Also, in the case
$(C,A)$ is the orbifold quotient of $\PP^1$, we may replace $(C,A)$
by $(\PP^1,\emptyset)$ without loss and then delete a point to form 
$(\CC,\emptyset)$.
Hence in all the nonhyperbolic cases, we may assume that $(C,A)$ is
uniformized by $\CC$, which we assume henceforth.

Let $r:\CC \rightarrow (C,A)$ be the uniformization map where the 
orbifold structure $A$ is defined by the minimum multiplicity on the fibers
of $f:X\rightarrow C$. Let $\tilde Y$ be the minimal resolution of 
singularities of the normalization of the fiber product of $Y$ with $\CC$
over $C$. $\tilde Y$ naturally fibers over $\CC$ as the pullback 
fibration of $Y$ over $C$ via $r$. 
Take $\tilde X$ to be the inverse image of $X$ in $\tilde Y$.
Then $\tilde X$ is a Zariski open subset of $\tilde Y$ with an
irreducible component of
multiplicity one at every fiber by construction. We will denote the 
union of all such irreducible components of  $\tilde X$ minus the singular
points of the fibers of $\tilde X$
by $\tilde X^\circ$, which is a Zariski open subset of $\tilde X$. 
It surjects to $\CC$ under the fibration map. We get the following
diagram of maps:
\begin{equation} \label{diag:ell}
\begin{array}{ccrcrcl}
\tilde X^\circ &  \hookrightarrow & \tilde X & \hookrightarrow 
& \tilde Y  & \rightarrow &  \CC \ \\
& & \tilde r \downarrow \, & &  \downarrow \, & & \,\downarrow   r \\
& & X & \hookrightarrow & Y  & \rightarrow &  C  
\end{array}
\end{equation}

By repeatedly blowing down the ($-1$)-curves in the fibers of $\tilde Y$, we
obtain a relatively minimal fibration $\tilde Y_{min}$ over $\CC$
(a fibration without such curves on its fibers). Moreover, 
proposition III.8.4 in \cite{BPV} (or rather its proof
plus the  fact that the arithmetic genus of a multiple $\PP^1$
fiber must be negative, which is disallowed, by II.11(b) of \cite{BPV})
tells us that $\tilde Y_{min}$ is
unique unless all of its fibers over $\CC$ are $\PP^1$, which implies
that $\tilde Y_{min}$ is a trivial $\PP^1$
bundle by the Grauert-Fischer theorem and hence also unique. Since the generic
fiber $X_c$  of $X$ is nonhyperbolic, $X_c$ is either $\PP^1$ with 
at most two punctures or it is an elliptic curve. In the former case
$\tilde Y_{min}$ is a trivial $\PP^1$ bundle over $\CC$
while in the latter $\tilde Y_{min}$ is a relatively minimal elliptic
fibration. We will call these cases case (1) and case (2), respectively,
and consider them separately below.
\medskip

The key is to note that, by construction, 
points of $\tilde X^\circ$ are infinitely near
points on $\tilde Y_{min}$ (including zeroth order) representing jets of
local sections of $\tilde Y_{min}$ over $\CC$. 
Thus, if we let $\tilde Y_{min}^\circ$ be the union of
all irreducible components of multiplicity one of fibers of $\tilde Y_{min}$
minus the singular
points of the fibers of $\tilde Y_{min}$, then  $\tilde Y_{min}^\circ$
surjects to $\CC$.
In particular, $\tilde Y_{min}^\circ$ does not have any multiple 
fiber and contains the image of $\tilde X^\circ$.
We have the following diagram:
$$ \tilde p: \tilde X\stackrel{\iota}{\hookrightarrow} \tilde Y
\stackrel{j}{\rightarrow} \tilde Y_{min} \stackrel{\pi}{\rightarrow} \CC$$
To conclude our proof,
it is sufficient to construct a dominating holomorphic map from 
$\CC^2$ to $\tilde X$ since a dominating holomorphic map 
to $X$ can  then be obtained by composing with the natural 
surjective holomorphic map $\tilde r$
from $\tilde X$ to $X$.
\medskip

\noindent
{\bf Case (1):} The generic fiber
$X_c$ is $\PP^1$ with at most two punctures: 
\smallskip

In the case of exactly two
punctures, $\tilde X$ is isomorphic via $j\circ\iota$
to the complement of a double section $D$ in 
the $\PP^1$ bundle $\tilde Y_{min}$  
over the complement of a discrete set $E$ in $\CC$. 
Since $\tilde p$ maps $\tilde X^\circ$
onto $\CC$, we may choose a point in $\tilde X^\circ_s$ for each $s\in E$,
thereby defining a jet prescription on $\tilde Y_{min}$ over $E$. Applying
theorem~\ref{thm:jet} then gives us a 
holomorphic map $F$ from $\CC^2$ to $\tilde Y_{min}$ that surjects to the
$\tilde Y_{min}\setminus D$ above $\CC\setminus E$. Moreover, $F$ is fiberwise
(meaning $\pi\circ F = p_1$, the projection to the first factor of $\CC^2$)
and maps sections of $p_1$ to sections of $\pi$ satisfying the 
jet prescription on $E$ given above. Hence $F$ lifts to a dominating 
holomorphic map
to $\tilde X$ which then gives a dominating map to $X$ by composition.
\smallskip

The case of no puncture can be reduced to that of a puncture 
by removing a section from the trivial $\PP^1$ 
bundle  $\tilde Y_{min}$. Hence 
it remains only to consider the case of exactly one puncture to establish
case (1). This puncture can be regarded as the infinity section $D_\infty$
in the trivial $\PP^1$ bundle  $\tilde Y_{min}$, whose complement
is a trivial line bundle $L$ over $\CC$ as $\CC$ is Stein. Again,
over the complement of a discrete set $E$ in $\CC$, 
$\tilde X$ is isomorphic via $j\circ\iota$ to this line bundle.

We now choose a point in $\tilde X^\circ_s$ for each $s\in E$,
thus defining a jet prescription over $E$. Using
a trivialization of  $\tilde Y_{min}$, this 
jet prescription corresponds to prescribing a Laurent polynomial at 
each point of $E$, and the basic interpolation theorems as before allow
us to find a meromorphic function $q(z)$ on $\CC$ satisfying this jet
prescription on $E$. Let $p(z)$ be an entire function vanishing to an
order greater than that of the jet prescription at each point of $E$
and nonzero outside $E$. If $w$ is the 
vertical coordinate of $\CC^2$, then the map
$$ (z,w)\mapsto (z, p(z)w+q(z)) $$
gives a holomorphic dominating map $F$ from $\CC^2$ to  
$\tilde Y_{min}=\CC\times \PP^1$ which commutes with the projection 
to the first factor (i.e., $F$ is a fiberwise map).
It certainly maps onto $L\hookrightarrow \tilde Y_{min}$ outside
$E$ and satisfies the jet prescriptions over $E$ for the image of all sections
in $\CC^2$. Hence $F$ lifts to a dominating holomorphic map to $\tilde X$
as required.

\medskip
\noindent
{\bf Case (2):} The generic fiber $X_c$ is an elliptic curve. Then, as
before, $\tilde X^\circ$ is isomorphic to $\tilde Y_{min}^\circ$
via $j\circ\iota$ above the complement of a discrete set $E$ in $\CC$, 
above which 
$\tilde Y_{min}^\circ$ is in particular a smooth elliptic fibration.
Choosing a point in  $\tilde X^\circ$ above each point in $E$ gives a
jet prescription of local sections through a component of  
$\tilde Y_{min}^\circ$ above each point in $E$.
If we remove all the other components of $\tilde Y_{min}^\circ$, then
we will show that
what remains is isomorphic to a Jacobian fibration, which is a 
holomorphic quotient of $\CC^2$ where
the quotient map is fiberwise and is a local isomorphism.
Assuming we have shown this for the moment, then our jet prescription in 
$\tilde Y_{min}^\circ$ above a point $s\in E$ gives rise to possibly
infinitely many jet prescriptions in $\CC^2$ above $s$. Choosing one
such jet prescription for each $s\in E$ and arguing as in 
the second part of case (1) above gives us a dominating map from $\CC^2$
to itself satisfying this jet prescription on $\CC^2$. Hence, composing
this dominating map with the quotient map above gives us a dominating map
to $\tilde Y_{min}^\circ$ satisfying the jet prescription given there.
This map would then lift to a dominating map to $\tilde X^\circ$, establishing
our theorem~\ref{thm:surf}.
\smallskip

So it remains to establish the isomorphism to a Jacobian
fibration as mentioned above. We will do this
in complete analogy with the proof of
theorem~\ref{thm:Jac} in section 3. Fortunately, Kodaira in \cite{Kod}
has worked out the analogue of our proposition~\ref{prop:L3} which we
quote as follows (see also \cite[\S V.9]{BPV}).

\begin{prop}   \label{prop:Jel}
Let $f : X \rightarrow C$ be a relatively minimal elliptic fibration
over a curve $C$ with a holomorphic
section $\sigma : \ C \rightarrow X$.  Let $X
_{\sigma} ^ {\prime}$ consists of all irreducible components of fibers
$X_s$ not meeting $\sigma (C)$, and let $X^{\sigma} = X \setminus X
_{\sigma} ^ {\prime}$.  Then there is a canonical fiber-preserving
isomorphism $h$ from {\rm Jac}$( f )$ onto $X^{\sigma}$ mapping the
zero-section in {\rm Jac}$( f )$ onto $\sigma (C)$.
\end{prop}

Since Jac$(f)$ is a holomorphic fiberwise quotient of the trivial line
bundle $ {L}: = {f} _{\ast 1} {\cal O} _{X}$ 
over $\CC$ (as in section 3, see also \cite{bl}), 
the proposition implies that as soon as $f$ has a section,
sections of this line bundle act on $X$ by translation in the fiber direction
via addition.

This proposition allows us to finish the proof by constructing a 
section $\sigma$ of $\tilde Y_{min}^\circ$ intersecting the component 
corresponding to our jet prescription for each $s\in E$, since
deleting the other components from $\tilde Y_{min}^\circ$ gives
us exactly $\tilde Y_{min}^\sigma$ by construction. This section 
$\sigma$ exists
by the following lemma. Its proof is verbatim that at the end of 
section 3, where we emulated Oka's original approach 
(\cite{Oka}) to the Oka-Grauert
principle for bundle of Lie groups, 
and is fitting as a conclusion to this section.

\begin{lemma} \label{lem:sigma} Given a relatively minimal elliptic fibration 
$f: X\ra C$ without multiple fibers, assume $C$ is non-compact.
Then $f$ has a holomorphic section. Furthermore, if we choose an 
irreducible component with multiplicity one in each fiber and 
let $X^\circ$ be their union minus the singular points of the fibers,
then the above holomorphic section can be chosen to have image in $X^\circ$.
\end{lemma}

\noindent{\bf Proof:} From  Kodaira's table of non-multiple
singular fibers (\cite{Kod} or \cite[Table 3 p. 150]{BPV}), we see that
every fiber which is not multiple in a relatively minimal elliptic fibration
has a component of multiplicity one.  We choose one such for each fiber 
of $f$ to construct  $X^\circ$ as in the lemma. Then, every point on
${C}$ admits a neighborhood with a section in $X^\circ$.
The rest of the proof is now verbatim that given at the end of
section 3 replacing proposition~\ref{prop:L3} by 
proposition~\ref{prop:Jel}.
$\BOX$

\section{Further remarks}
Our main theorems on dominating maps
allow us to show for example that blowing up
$\PP^2$ at $x$ and taking the complement of the strict transform of the 
normal crossing quartic consisting of two lines intersecting at $x$ and
a conic gives us an algebraic surface that is 
dominated holomorphically by $\CC^2$. This despite the fact that the 
complement of the same quartic in $\PP^2$ is not so dominable, as
shown in \cite{CG} using a well known lemma of Kodaira (\cite{K}).
\smallskip

As for  theorem~\ref{thm:section} on the existence of a section,
we first note that it and hence also theorem~\ref{thm:dominate} 
are false for a compact base
$R$ since the double section $D$ can be ample. It is not difficult by our 
method to characterize those double sections on a 
compact ruled surface which admit 
sections disjoint from it. The same characterization when $D$ is replaced by
a section is classical and can be found,
for example, in \cite[\S V.2]{Ha}.

We can also consider the possibility of other generalizations of
theorem~\ref{thm:section}, especially various versions of the
``$n$-section'' problem. However, the following problem already has a
negative answer for $n=3$. 

\begin{question}
Let $D$ be an $n$-section of a $\PP^1$ bundle over a noncompact Riemann
surface. That is, $D$ is a reduced divisor intersecting each fiber exactly
$n$ times counting multiplicity.
Does there exist a section of the bundle avoiding $D?$
\end{question}

We have already shown that the answer is yes for $n=2$. 
To see that the answer is no for $n=3$ , we consider the following example:
\smallskip

Let $L_0(z) = 0$, $L_1(z) = z$, and let $T$ be the triple section in $\CC
\times \PP^1$ given by the union of the infinity section and the graphs of
$L_0$ and $L_1$ ($L_0$ and $L_1$ could be replaced by any two lines meeting
at the origin).  Suppose $\sigma$ is an entire function whose graph avoids
$T$.  We will exclude the fiber over $0 \in \CC$ and pull-back with the
covering map to get a contradiction. 
\smallskip

Let $\Phi(u,v) = (e^u, v)$, and let $T^*$ be $\Phi^{-1}(T)$.  Then $T^*$ 
is the union of the infinity section and the graphs of $L_0^*$ and 
$L_1^*$, where $L_0^*(u) = 0$ and $L_1^*(u) = e^u$.  Also, $\Phi^{-1}$ of 
the section $\sigma$ is the graph of $\sigma^*(u) = \sigma(e^u)$, and 
this graph is contained in the complement of $T^*$.  

Finally, let $\Psi(u,v) = (u, v/e^u)$.  Then $\Psi(T^*)$ is $\CC \times 
\{0, 1, \infty\}$, and hence $\Psi$ applied to the graph of $\sigma^*$ is 
the graph of an entire function $\sigma^{**}$ which never attains the
values 0 or 1 and which hence must be 
some constant $c \neq 0,1$.  But also $\sigma^{**}(u) = \sigma^*(u)/e^u = 
\sigma(e^u)/e^u$.  

Thus $\sigma(e^u) = ce^u$, or equivalently $\sigma(z) = cz$.  But then
$\sigma(0) =  0$, which contradicts the fact that the graph of
$\sigma$ avoids $T$. 
\smallskip

Note that this is an example over a nonhyperbolic base. The same question
for a hyperbolic base curve remains an important open problem.
\bigskip

As for nonhyperbolically fibered surfaces, it is clear that dominability
by $\CC^2$ is a bimeromorphic invariant. It follows that the hyperbolicity
of the orbifold base curve from the minimal multiplicities
of components of the fibration is also invariant, as is the property of having
a Zariski dense holomorphic image of $\CC$ by theorem~\ref{thm:surf}. 
However, it is still a 
problem to find the connection of this invariant 
to other natural bimeromorphic invariants of the surface, especially
to tensorial invariants, as one would like to find a differential
geometric connection.

In this connection, \cite{bl} gives a 
fundamental group characterization of dominability for 
projective algebraic surfaces modulo certain nonelliptic K3 surfaces,
which we quote below.

\begin{theorem}  \label{thm:2}
A projective surface $X$ not birationally equivalent to a K3 surface
is dominable by $\CC^2$ if and only if 
it has Kodaira dimension less than two and its fundamental group
is a finite extension of an abelian group (of even rank four or less).
If $\kappa (X)=-\infty$, then the fundamental group condition can be
replaced by the simpler condition of non-existence of more than one 
linearly independent holomorphic one-form. If $\kappa (X)=0$ and $X$
is not birationally equivalent to a K3 surface, then $X$ is dominable
by $\CC^2$. If $X$ is birationally equivalent to an elliptic K3
surface or to a Kummer K3 surface, then $X$ is dominable by $\CC^2$.
\end{theorem}

The proof of this characterization 
depends on the fact that in the case of a properly
fibered surface, every nonmultiple fiber has an irreducible component
of multiplicity one. This is no longer true when the fibers 
of a fibered surface are allowed to be noncompact. In fact, the following
is a simple example of a 
fibered, simply connected, quasiprojective surface that is
generically ruled (a $\PP^1$ bundle outside some finite subset) but is
not dominable by $\CC^2$.
\smallskip

Consider the fibration $\CC\times\PP^1\rightarrow \CC$
obtained by projecting to the first factor.  Note that the base curve $C$
is quasiprojective. Consider
two points $a=(0,2)$, $b=(0,3)$ in the fiber at $0$ and the following 
operation on this surface. Take two consecutive blowups at $a$ 
(first blow up at $a$ and then blow up at the singular point of
the resulting fiber) and three
such blowups at $b$ to obtain an irreducible component of multiplicity
$2$ above $a$ and one of multiplicity $3$ above $b$ for the resulting fiber.
Delete all other components of this fiber from the blown-up surface
and call the resulting open surface $Y^*$.

It is easy to see that the resulting surface is simply connected: 
First note that any loop can be homotoped to the complement of the fiber
(since the fiber is real codimension two in the surface)
and so corresponds to a multiple of the generator $\alpha$ of 
$\pi_1 (\CC\setminus \{0\})$.  Take
a simple loop $S$ outside the fiber corresponding to  $\alpha$.
Compose $S$, with a base point fixed, with two loops outside 
the fiber, one corresponds to $-3\alpha$ in a disk
through $b$ transversal to the fiber and the other corresponds
to $2\alpha$ in a disk through $a$ transversal to the fiber. 
Since these two loops are  contractible in
$Y^*$, the composition is homotopic to $S$ but corresponds to
the zero multiple of $\alpha$ and hence is homotopically trivial as the
generic fiber is simply connected. So, $Y^*$ is simply connected. 
Note that the fiber $Y^*_0$ is not 
a multiple fiber as the multiplicities of its components are relatively 
prime.

Now we do exactly the same operation to two additional fibers giving
a surface $Y^{**}$ which is simply connected by the same argument
as above. However, the orbifold Euler characteristic for the base
orbifold curve using the minimal multiplicity of the components of the
fibers of  $Y^{**}$ is $1-3(1-1/2)=-1/2<0$ showing that the orbifold
curve is hyperbolic. Hence any holomorphic image of $\CC$ 
in  $Y^{**}$ must lie in a fiber as a consequence of our Ahlfors'
lemma, theorem~\ref{Ahlfors}.


\noindent
Gregery T. Buzzard\\
Department of Mathematics\\
Cornell University\\
Ithaca, NY  14853
USA

\medskip
\noindent
Steven Shin-Yi Lu\\
Department of Mathematics\\
University of Waterloo\\
Waterloo, Ontario, N2L3G1\\
Canada

\end{document}